\documentclass{article}

\usepackage{amsmath}
\usepackage{amssymb}
\usepackage{epsfig}
\usepackage{tabularx}

\newcommand{\x}{\textbf{x}}
\renewcommand{\v}{\textbf{v}}
\renewcommand{\k}{\textbf{k}}
\renewcommand{\d}{\partial}
\newcommand{\R}{\mathbb{R}}
\renewcommand{\u}{\textbf{u}}
\newcommand{\eps}{\varepsilon}
\newcommand{\sinc}{\text{sinc}}
\newcommand{\Qt}{\tilde{Q}}
\newcommand{\dt}{\Delta t}
\newcommand{\dx}{\Delta x}

\begin{document}

\title{High Performance Computing with a Conservative Spectral Boltzmann Solver \thanks{This paper has been accepted to the Proceedings of the 28th International Symposium on Rarefied Gas Dynamics}}


\author{Jeffrey R. Haack \thanks{Department of Mathematics, The University of Texas at Austin, 2515 Speedway, Stop C1200 Austin, Texas 78712}  \and Irene M. Gamba \thanks{Department of Mathematics, The University of Texas at Austin, 2515 Speedway, Stop C1200 Austin, Texas 78712 and ICES, The University of Texas at Austin, 201 E. 24th St., Stop C0200, Austin, TX 78712} 
}

\maketitle

\begin{abstract}
We present new results building on the conservative deterministic spectral method for the space inhomogeneous Boltzmann equation developed by Gamba and Tharkabhushaman. This approach is a two-step process that acts on the weak form of the Boltzmann equation, and uses the machinery of the Fourier transform to reformulate the collisional integral into a weighted convolution in Fourier space. A constrained optimization problem is solved to preserve the mass, momentum, and energy of the resulting distribution. We extend this method to second order accuracy in space and time, and explore how to leverage the structure of the collisional formulation for high performance computing environments. The locality in space of the collisional term provides a straightforward memory decomposition, and we perform some initial scaling tests on high performance computing resources. We also use the improved computational power of this method to investigate a boundary-layer generated shock problem that cannot be described by classical hydrodynamics. 
\end{abstract}


\section{Introduction} \label{sec:intro}

There are many difficulties associated with numerically solving the
Boltzmann equation, most notably the dimensionality of the problem and
the conservation of the collision invariants. For physically relevant
three dimensional applications the distribution function is seven
dimensional and the velocity domain is unbounded. In addition, the
collision operator is nonlinear and requires evaluation of a five
dimensional integral at each point in phase space. The collision
operator also locally conserves mass, momentum, and energy, and any
approximation must maintains this property to ensure that macroscopic
quantities evolve correctly.

Spectral methods are a deterministic approach that compute the
collision operator to high accuracy by exploiting its Fourier
structure. These methods grew from the analytical works of Bobylev \cite{Bob88} developed for the Boltzmann equation for 
 Maxwell type potential interactions and integrable angular cross section, where the corresponding Fourier transformed equation takes a 
closed form. Spectral approximations for this type of models where first proposed by Pareschi and Perthame \cite{ParPer96}.
Later Pareschi and Russo \cite{ParRus00} applied this work to variable hard potentials by periodizing the problem and its solution and implementing spectral collocation methods.

These methods require $O(N^{2d})$ operations per evaluation of the collision operator, where $N$ is the
total number of velocity grid points in each dimension. While convolutions 
can generally be computed in $O(N^d \log N)$ operations, the
presence of the convolution weights requires the full $O(N^{2d})$ computation of the convolution ,  except for a few special cases, e.g.,
the Fokker-Planck-Landau collision operator \cite{ParRusTos00, MouPar06}. Spectral methods provide many advantages over Direct Simulation Monte Carlo Methods (DSMC) because they are more suited to time dependent problems, low Mach number flows, high mean velocity flows, and flows that are away from equilibrium. In addition, deterministic methods avoid the statistical fluctuations that are typical of particle based methods.

Inspired by the work of Ibragimov and Rjasanow \cite{IbrRja02}, Gamba
and Tharkabhushanam \cite{GamTha09, GamTha10} observed that the
Fourier transformed collision operator takes a simple form of a weighted
convolution and  developed a spectral method based
on the weak form of the Boltzmann equation that provides a general
framework for computing both elastic and inelastic collisions. Macroscopic conservation is enforced by solving a numerical constrained optimization problem that finds the closest distribution function in $L_2$ to the output of the collision term that conserves the macroscopic quantities.  These methods do not rely on periodization but rather on the use of the FFT tool in the computational domain, where convergence to the solution of the continuous problem is obtained by the use of the extension operator in Sobolev spaces.

This paper presents extensions to this method and its implementation on high performance computing resources. Without loss of generalization, we restrict this presentation to elastic collisions. We present a second order in time and space extension of the Gamba and Tharkabhushanam method allowing for nonuniform grids. This method has been implemented on the Lonestar supercomputer at the Texas Advanced Computing Center (TACC) and timing studies are provided to explore the scaling of the method to more and more processors for large problems. Finally we present some $1D$ in physical space results for a problem proposed by Aoki et al.\cite{AokiSoneNisSug91} where a sudden change in wall temperature results in a shock that cannot be explained by classical hydrodynamics.

\section{The space inhomogeneous Boltzmann equation} \label{sec:BTE}
The space inhomogeneous initial-boundary value problem for the Boltzmann equation is given by
\begin{equation}\label{BTE}
\frac{d}{dt} f(\x,\v,t) + \v \cdot \nabla_x f(\x, \v, t) = \frac{1}{\eps} Q(f,f),
\end{equation}
with
\begin{align*}
&\x \in \Omega \subset \R^d, \qquad \v \in \R^d \\
&f(\x,\v,0) = f_0(\x,\v) \\
&f(\x,\v,t) = f_B(\x,\v,t), \qquad \x \ \in \d \Omega.
\end{align*}
where $f(\v,t)$ is a probability density distribution in $\v$-space and $f_0$ is assumed to be locally integrable with respect to $\v$ and the spatial boundary condition $f_B$ will be specified in below. The dimensionless  parameter $\eps > 0$ is the scaled Knudsen number, which is defined as the ratio between the mean free path between collisions and a reference macroscopic length scale.
 
The collision operator $Q(f,f)(\x,\v,t)$ is a bilinear integral form local in $t$ and $\x$ and is given by
\begin{equation}\label{Q_general}
Q(f,f)(\cdot,\v,\cdot) = \int_{\v_\ast \in \R^d} \int_{\sigma \in S^{d-1}} B(|\v - \v_\ast|,\cos \theta) (f(\v_\ast ')f(\v') - f(\v_\ast)f(\v)) d\sigma d\v_\ast,
\end{equation}
where the velocities $\v', \v_\ast'$ are determined through a collision rule \eqref{velocity_interact}, depending on $\v, \v_\ast,$, and the positive term of  the integral in \eqref{Q_general} evaluates  $f$ in the pre-collisional velocities that will take the direction $\v$ after an interaction.  The collision kernel $B(|\v - \v_\ast|,\cos \theta)$ is a given non-negative function depending on the size of the relative velocity $\u := \v - \v_\ast$ and $\cos \theta = \frac{\u \cdot \sigma}{|\u|}$, where $\sigma$ in the $n-1$ dimensional sphere $S^{n-1}$ is referred as the scattering direction  
of the post-collisional elastic relative velocity.

For the following we will use the velocity elastic (or reversible) interaction law in center of mass-relative velocity coordinates
\begin{align} 
&\v' = \v + \frac{1}{2}(|\u|\sigma - \u), \qquad \v_\ast ' = \v_\ast - \frac{1}{2}(|\u|\sigma - \u) \label{velocity_interact} \\
&B(|\u|,\cos \theta) = |\u|^\lambda b(\cos \theta) \, \nonumber .
\end{align}
We assume that the differential cross section function $b(\cos \theta)$ is integrable with respect to $\sigma$ on $S^{d-1}$, referred to as the Grad cut-off assumption, and that it is renormalized such that 
\begin{equation} \label{cut-off}
\int_{S^{d-1}} b(\cos \theta) d\sigma = 1.
\end{equation}
The parameter $\lambda$ regulates the collision frequency as a function of the relative speed $|\u|$. This corresponds to the interparticle potentials used in the derivation of the collisional kernel and are referred to as variable hard potentials (VHP) for $0 < \lambda < 1$, hard spheres (HS) for $\lambda = 1$, Maxwell molecules (MM) for $\lambda = 0$, and variable soft potentials (VSP) for $-3 < \lambda < 0$. The $\lambda = -3$ case corresponds to a Coulombic interaction potential between particles. If $b(\cos \theta)$ is independent of $\sigma$ we call the interactions isotropic (like the case of hard spheres in three dimensions). 
%
%

\subsection{Boundary conditions} \label{sec:BC}
On the spatial boundary $\d \Omega$ we use a diffusive Maxwell boundary condition which is given by, for $\x \in \d \Omega$,
\begin{align} \label{BC}
f(\x,\v,t)  &= \frac{\sigma_w}{(2\pi R T_w)^(d/2)} \exp \left( - \frac{|\v - \textbf{V}_w|^2}{2 R T_w}\right), \qquad (\v - \textbf{V}_w)\cdot \textbf{n} > 0 \\
\sigma_w &= - \left(\frac{2\pi}{R T_w}\right)^{1/2} \int_{(\v - \textbf{V}_w)\cdot\textbf{n} < 0} (\v - \textbf{V}_w)\cdot\textbf{n} f(\x,\v,t) d \v, \nonumber
\end{align}
where $\textbf{V}_w$ and $T_w$ are the wall velocity and temperature, respectively, and $\textbf{n}$ is the unit normal vector to the boundary, directed into $\Omega$. The term $\sigma_w$ accounts for the amount of particles leaving the domain and ensures mass conservation in $\Omega$.


\subsection{Spectral formulation} \label{sec:spectral_cont}
The key step our formulation of the spectral numerical method is the use of the weak form of the Boltzmann collision operator. For a suitably smooth test function $\phi(\v)$ the weak form of the collision integral is given by
\begin{equation} \label{collision_weakform}
\int_{\R^d} Q(f,f) \phi(\v) d\v = \int_{\R^d \times \R^d \times S^{d-1}} f(\v)f(\v_\ast) B(|\u|,\cos \theta) (\phi(\v') - \phi(\v)) d\sigma d\v_\ast d\v
\end{equation}
If one chooses 
\begin{equation*}
\phi(\v) = e^{-i \zeta \cdot \v} / (\sqrt{2\pi})^d,
\end{equation*}
then \eqref{collision_weakform} is the Fourier transform of the collision integral with Fourier variable $\zeta$:
\begin{align}
\widehat{Q}(\zeta) &= \frac{1}{(\sqrt{2\pi})^d} \int_{\R^d} Q(f,f) e^{-i \zeta \cdot \v} d\v \nonumber \\
&= \int_{\R^d \times \R^d \times S^{d-1}} f(\v)f(\v_\ast) \frac{B(|\u|,\cos \theta)}{(\sqrt{2\pi})^d} (e^{-i \zeta \cdot \v'} - e^{-i \zeta \cdot \v}) d\sigma d\v_\ast d\v\\
&= \int_{\R^d} G(\u,\zeta) \mathcal{F}[f(\v)f(\v-\u)](\zeta) d\u,
\end{align}
where $\widehat{[\cdot]} = \mathcal{F}(\cdot)$ denotes the Fourier transform and 
\begin{equation} \label{G_eqn}
G(\u,\zeta) = |\u|^\lambda \int_{S^{d-1}} b(\cos \theta) \left(e^{-i\frac{\beta}{2} \zeta \cdot |\u|\sigma}e^{i\frac{\beta}{2} \zeta \cdot \u} - 1\right) d\sigma
\end{equation}
Further simplification can be made by writing the Fourier transform inside the integral as a convolution of Fourier transforms:
\begin{align} \label{Cont_spectral}
\widehat{Q}(\zeta) 
&= \int_{\R^d} \widehat{G}(\xi,\zeta) \hat{f}(\zeta - \xi) \hat{f}(\xi) d\xi,
\end{align}
where the convolution weights $\widehat{G}(\xi,\zeta)$ are given by
\begin{align} \label{Ghat_eqn}
\widehat{G}(\xi,\zeta) &= \frac{1}{(\sqrt{2\pi})^d}  \int_{\R^d} G(\u,\zeta) e^{-i \xi \cdot u} d\u 
\end{align}
These convolution weights can be precomputed once to high accuracy and stored for future use. For many collision types the complexity of the integrals in the weight functions can be reduced dramatically through analytical techniques. In this paper we will only consider isotropic scattering in dimension 3  $(b(\cos \theta) = 1/ 4\pi)$. In this case we have that 


\begin{align} \label{GHat_calc}
\widehat{G}(\xi,\zeta) =  \frac{4\pi}{(\sqrt{2\pi})^3} \int_{\R^+} r^{\lambda + 2} \left( \sinc\left(\frac{\beta r |\zeta|}{2}\right) \sinc(r|\xi - \frac{\beta}{2}\zeta|) - \sinc(r|\xi|)\right) dr.
\end{align}

This integral will be cut off at a point $r = r_0$, which will be determined below. Given this cutoff point, we can explicitly compute $\widehat{G}$ for integer values of $\lambda$. 
For other values of $\lambda$, this is simply a one-dimensional integral that can be precomputed to high accuracy using numerical quadrature. The entirety of the collisional model being used is encoded in the weights, which gives the algorithm a large degree of flexibility in implementing different models. 


\section{The Conservative Numerical Method} \label{sec:numerics_setup}

\subsection{Temporal and velocity space discretization} \label{sec:discretization}

We use an operator splitting method to separate the mechanisms of collisions and advection. The system is split into the subproblems
\begin{align}
\frac{\d}{\d t} f + v \cdot \nabla_\x f = 0 \\
\frac{\d}{\d t} f = Q(f,f),
\end{align}
which are solved separately. 

Each system is evolved in time using a second-order Runge-Kutta method, and the systems are combined using Strang splitting.

We assume that the distribution function is negligible outside of a ball 
\begin{equation} \label{Ball_domain_v}
B_{R_x}(\textbf{V}(\x)) = \{\v \in \R^d : |\v - \textbf{V(\x)}| \le R_x \},
\end{equation}
where $\textbf{V}(\x)$ is the local flow velocity which depends in the spatial variable $\x$. For ease of notation in the following we will work with a ball centered at $0$ and choose a length $R$ large enough that $B_{R_x}(\textbf{V}(\x)) \subset B_R(0)$ for all $\x$.

With this assumed support for the distribution $f$, the integrals in \eqref{Cont_spectral} will only be nonzero for $\u \in B_{2R}(0)$. Therefore, we set $L=2R$ and define the cube
\begin{equation} \label{Cube_domain_v}
C_L = \{ \v \in \R^d : |v_j| \le L,\,\, j = 1,\dots,d\}
\end{equation}
to be the domain of computation. For such domain, the computation of the weight function integral \eqref{GHat_calc} is cut off at $r_0=L$.

Let $N \in \mathbb{N}$ be the number of points in velocity space in each dimension. Then the uniform velocity mesh size is $\Delta v = \frac{2L}{N}$ and due to the formulation of the discrete Fourier transform the corresponding Fourier space mesh size is given by $\Delta \zeta = \frac{\pi}{L}$. 

The mesh points are defined by
\begin{align} \label{meshpoints}
v_\k &= \Delta v (\k - N/2) \nonumber \\
\zeta_\k &= \Delta \zeta (\k - N/2) \\
& \k = (k_1,\dots,k_d) \in \mathbb{Z}^d,\qquad 0 \le k_j \le N-1,\,\,\, j = 1,\dots,d
\end{align}


\subsection{Collision step discretization} \label{sec:collision_disc}
Returning to the spectral formulation \eqref{Cont_spectral}, the weighted convolution integral then becomes an integral over $-\frac{\pi}{\Delta v} \le \xi_j \le \frac{\pi}{\Delta v}, \, j = 1,\dots,d$. 
%

To simplify notation we will use one index to denote multidimensional sums with respect to an index vector $\textbf{m}$
\begin{equation*}
\sum_{\textbf{m}=0}^{N-1} = \sum_{m_1,\dots,m_d = 0}^{N-1}.
\end{equation*}

To compute $\widehat{Q}(\zeta_\k)$, we first compute the Fourier transform integral giving $\hat{f}(\zeta_k)$ via the FFT. The weighted convolution integral is approximated using the trapezoidal rule
\begin{equation}
\widehat{Q}(\zeta_\k)= \sum_{\textbf{m} = 0}^{N-1} \widehat{G}(\xi_{\textbf{m}},\zeta_\k) \hat{f}(\xi_{\textbf{m}}) \hat{f}(\zeta_\k - \xi_{\textbf{m}}) \omega_{\textbf{m}},
\end{equation}
where $\omega_\textbf{m}$ is the quadrature weight and we set $\hat{f}(\zeta_\k - \xi_{\textbf{m}}) = 0$ if $(\zeta_\k - \xi_{\textbf{m}})$ is outside of the domain of integration. We then use the inverse FFT on $\widehat{Q}$ to calculate the integral returning the result to velocity space. 

Note that in this formulation the distribution function is not periodized, as is done in the collocation approach of Pareschi and Russo \cite{ParRus00}. This is reflected in the omission of Fourier terms outside of the Fourier domain. All integrals are computed directly only using the FFT as a tool for faster computation.The convolution integral is accurate to at least the order of the quadrature. The calculations below use the trapezoid rule, but in principle Simpson's rule or some other uniform grid quadrature can be used. However, it is known that the trapezoid rule is spectrally accurate for periodic functions on periodic domains (which is the basis of spectral accuracy for the FFT), and the same arguments can apply to functions with sufficient decay at the integration boundaries \cite{Atkinson}. These accuracy considerations will be investigated in future work. The overall cost of this step is $O(N^{2d})$. 

\subsection{Discrete conservation enforcement} \label{sec:conservation}
This implementation of the collision mechanism does not conserve all of the quantities of the collision operator. To correct this fact, we formulate these conservation properties as a constrained optimization problem as proposed in \cite{GamTha09, GamTha10}. Depending on the type of collisions we can change this constraint set (for example, inelastic collisions do not preserve energy). We focus here just on the case of elastic collisions, which preserve mass, momentum, and energy. 

Let $M = N^d$ be the total number of grid points, let $\tilde{\textbf{Q}} = (\Qt_1, \dots, \Qt_M) ^T$ be the result of the spectral formulation from the previous section, written in vector form, and let $\omega_j$ be the quadrature weights over the domain in this ordering. Define the integration matrix
\begin{equation*}
\textbf{C}_{5\times M} = \left(\begin{array}{c} \omega_j \\ v_j^i \omega_j \\ |\v_j|^2 \omega_j \end{array} \right),
\end{equation*}
where $v^i,\, i=1,2,3$ refers to the $i$th component of the velocity vector. Using this notation, the conservation method can be written as a constrained optimization problem. 

\begin{equation} 
\text{Find } \textbf{Q} = (Q_1,\dots,Q_M)^T \text{ that minimizes } \frac12 \|\tilde{\textbf{Q}} - \textbf{Q}\|_2^2 \text{ such that } \textbf{C} \textbf{Q} = \textbf{0}
\end{equation}

The solution is given by
\begin{align}
\textbf{Q} &= \tilde{\textbf{Q}} + \textbf{C}(\textbf{C}\textbf{C}^T)^{-1} \textbf{C} \tilde{\textbf{Q}} \nonumber \\
&:= \textbf{P}_N \tilde{\textbf{Q}}
\end{align}

%
The overall cost of the conservation portion of the algorithm is a $O(N^d)$ matrix-vector multiply, significantly less than the computation of the weighted convolution.


\subsection{Spatial and Transport discretization} \label{sec:transport}

For simplicity this will be presented in 1D in space, though the ideas apply to higher dimensions. In this case the transport equation reduces to
\begin{equation*}
\frac{\d f}{\d t} (x,\v, t) + v_1 \frac{\d}{\d x} f(x,\v,t) = 0.
\end{equation*}

We partition the domain into cells of size $\dx_j$ (not necessarily uniform) with cell centers $x_j$. Using a finite volume approach, we integrate the transport equation over a single cell to obtain
\begin{align*}
\frac{f^{n+1}_j (v_i) - f^n_j (v_i)}{\dt} + \frac{F^n_{j+1/2} - F^n_{j-1/2}}{\dx_j} = 0,
\end{align*}
where $t^n = n\dt$ and $F^n_{j\pm1/2}$ is an approximation of the edge fluxes $v_1f$ of the cell between time $t^n$ and $t^{n+1}$. We use a second order upwind scheme defined by 
\begin{equation}\label{upwind}
F^n_{j+1/2} = \begin{cases} v_1 (f^n_j + \frac{\dx}{2}\sigma^n_j), \qquad &v_1 \ge 0 \\
		                         v_1(f^n_{j+1} - \frac{\dx}{2}\sigma^n_{j+1}), \qquad &\text{otherwise,} \end{cases}
\end{equation}
where $\sigma_j$ is a cell slope term used in the reconstruction defined by the minmod limiter.

On wall boundaries the incoming flux is determined using ghost cells and the diffusive reflection formula \eqref{BC}. For problems without meaningful boundary interactions (e.g. shocks), a no-flux boundary condition is applied for the incoming characteristics.

\section{Parallelization}
The major bottleneck in parallelizing a program is memory access times related to communication between processes. However, the relative locality of the dynamics of the Boltzmann equation allow for a straightforward decomposition of phase space. In each time step, a single grid point only ``sees"  the particles at the same spatial grid point, through the collision term, and particles with the same velocity at neighboring spatial grid points, through the transport term. As most of the computational time is spent on the collision term we choose to keep all of the information needed for this step on the same computational node, and thus partition the spatial grid points across the computational nodes. 

Further parallelization can be realized on each node by noting that the computation of $\widehat{Q}(\zeta)$ is simply a sum involving the precomputed weights $\widehat{G}$, the distribution function $\hat{f}$, and $\zeta$, and is completely independent of computation on another velocity grid point. We use OpenMP \cite{openMP08} to distribute the weighted convolutions to each computational core on the node. Each core computes the convolution for a single grid point in $\zeta$, marching through until all of the convolutions on the node are complete. Because there is no memory transfer required between cores, this should speed up the computation of the collisional terms by a factor of $p$, where $p$ is the total number of cores used in the computation.

We use the Message Passing Interface (MPI) \cite{MPI} to communicate distribution function data between nodes using an interleaved ghost cells technique \cite{GroppMPI}. Each node maintains two spatial ghost cells on the left and right sides of its local spatial domain, which are filled by the distribution functions sent from the neighboring node. Once this information is received, the node has enough information to update the ``regular" spatial grid points through the finite volume transport scheme described above.

To make a rough estimate of the total speedup let $n$ be the number of nodes used, and $np$ be the number of cores used (assume a fixed number of cores per node). Then the speedup can be described by, to leading order,
\begin{equation}
\frac{T_{\textrm{SERIAL}}}{T_{\textrm{PARALLEL}}} = \frac{ C N_x N^6 T_{\textrm{FLOP}} }{4nN^3 T_{\textrm{MEM}} + C N_x N^6 T_{\textrm{FLOP}}/np}, 
\end{equation}
where $T_{\textrm{MEM}}$ is the time to transfer a single floating point piece of data, and $T_{\textrm{FLOP}}$ is the time for a single floating point operation. As $n$ becomes large with $N$ and $N_x$ fixed, more and more data transfer is required, however one would need $4n^2 p T_{MEM}/T_{FLOP} \approx C N_x N^3$ before the memory access time would begin to dominate the collisional computations.


\section{Numerical results}\label{sec:numerics}
We test this method with the sudden change in wall temperature example suggested by Aoki et al. \cite{AokiSoneNisSug91}. In this problem the gas is initially at equilibrium and the temperature of the wall at the boundary of the domain is instantaneously changed at $t=0$. This gives rise to a discontinuous distribution function at the wall, which propagates into the domain and eventually forms a shock on the interior of the domain. This problem is especially well suited to a deterministic method because the discontinuity near the wall results in a distribution function that is far from equilibrium, which is difficult to simulate with Monte Carlo based solvers.

\begin{figure}[!htbp]\label{SuddenHeat}
\includegraphics[width=0.45\linewidth]{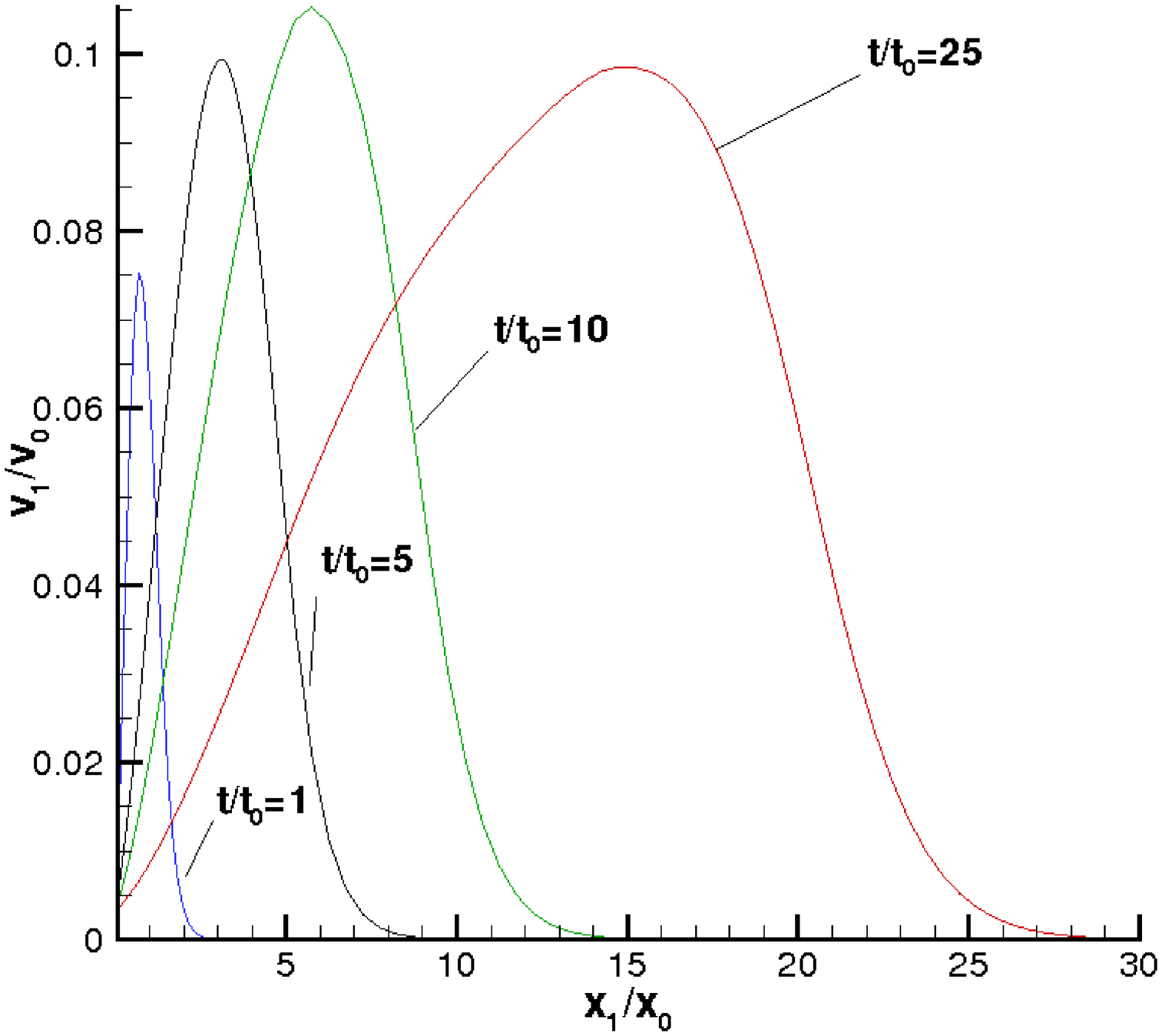} 
\includegraphics[width=0.45\linewidth]{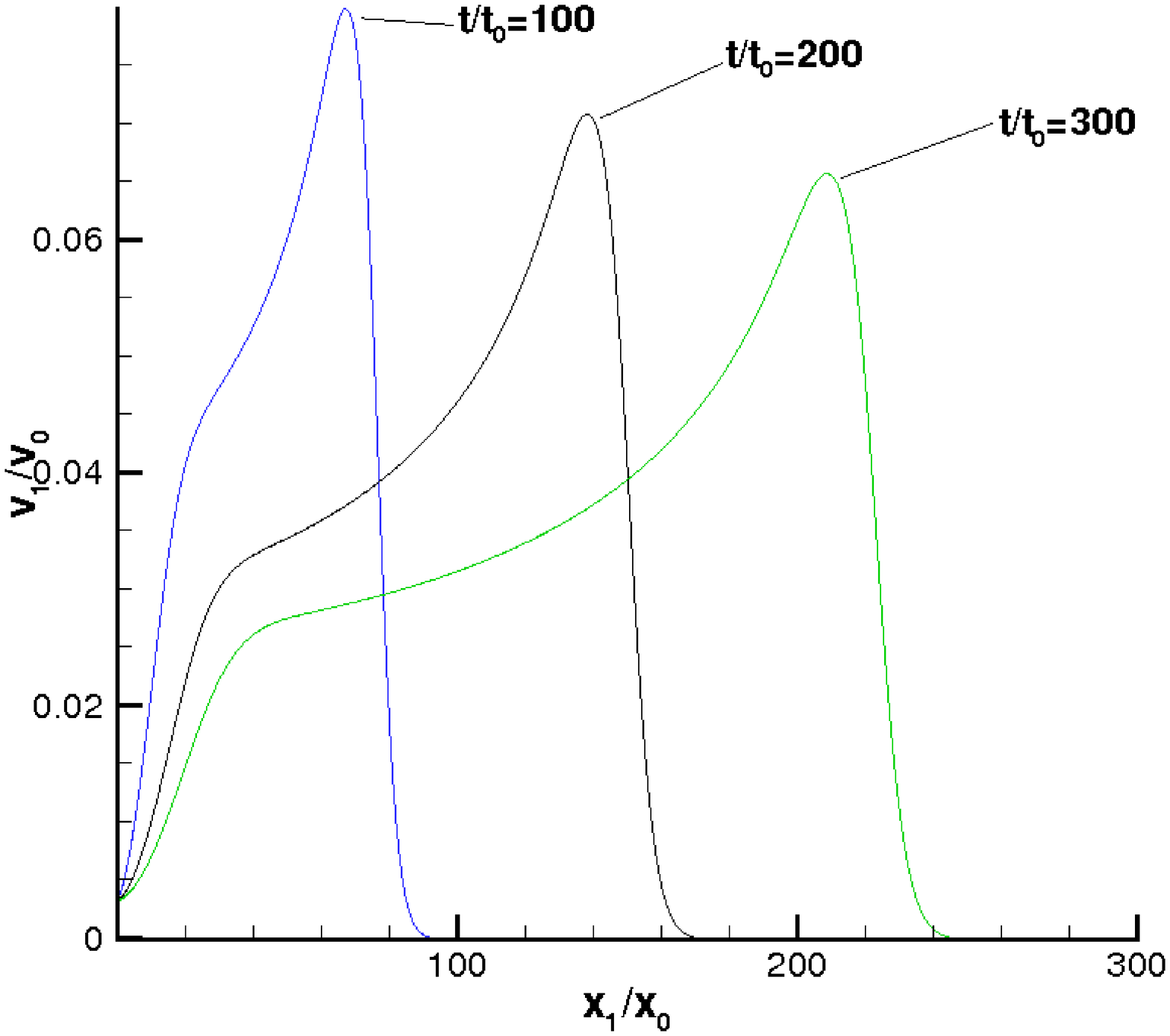} 
\caption{Formation of shock from suddent change in wall temperature.  Evolution of bulk velocity.}
\end{figure}

In Figure \ref{SuddenHeat} we show shock formation due to a sudden change in wall temperature. Unlike the computations in \cite{AokiSoneNisSug91}, only two grid points are used per mean free path in the interior of the domain. Near the wall, the grid is refined to eight points per mean free path to better capture the finer dynamics where the distribution function is discontinuous. In both examples, the scaled Knudsen number $\eps = 1$.  Note that despite the discontinuity we do not observe any Gibbs phenomena in the solution. We hypothesize that this is due to the mixing effects of the convolution weights; this will be explored in future work.

\begin{figure}[!htbp]\label{Marginal}
\includegraphics[width=0.45\linewidth]{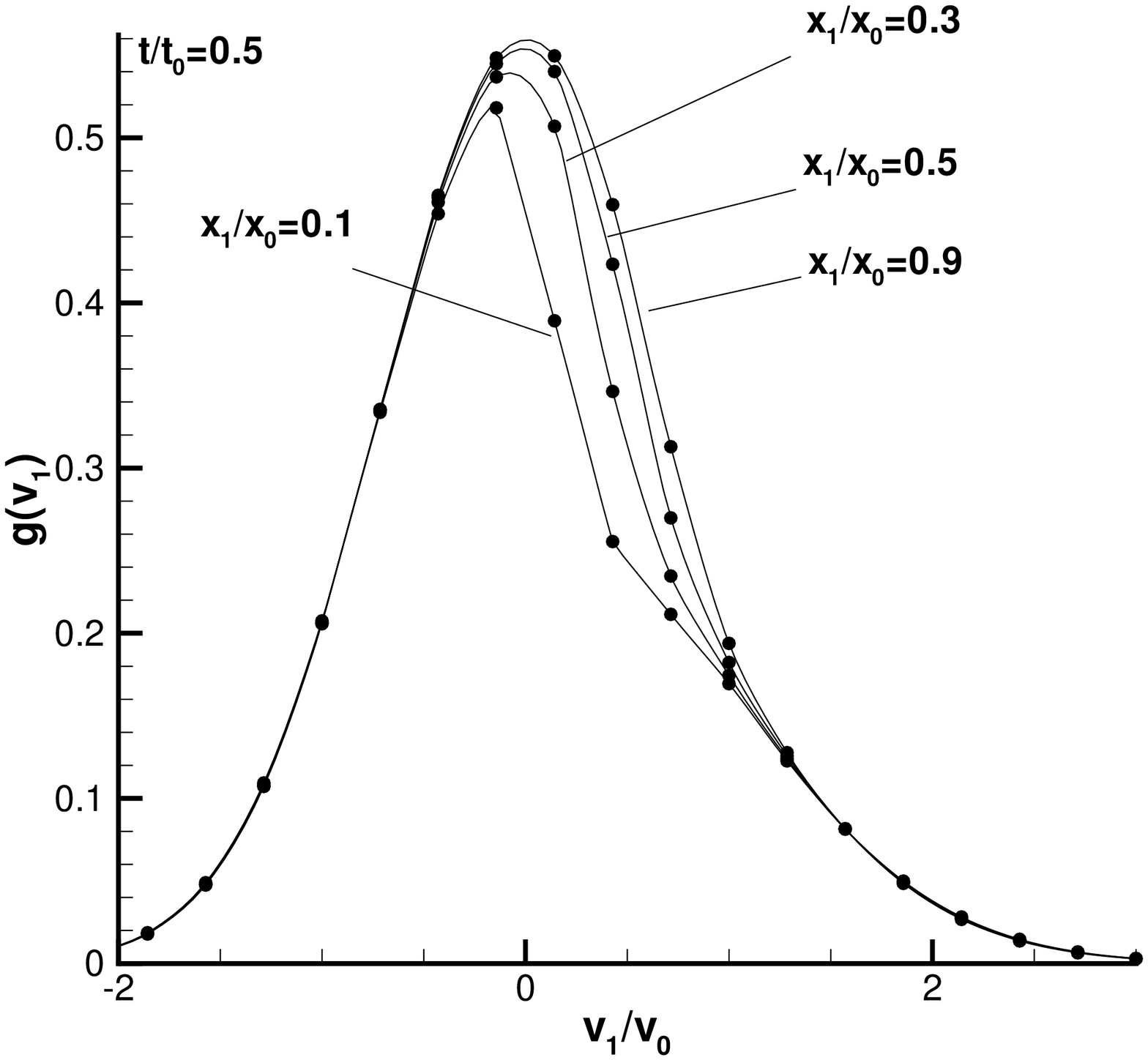} 
\includegraphics[width=0.45\linewidth]{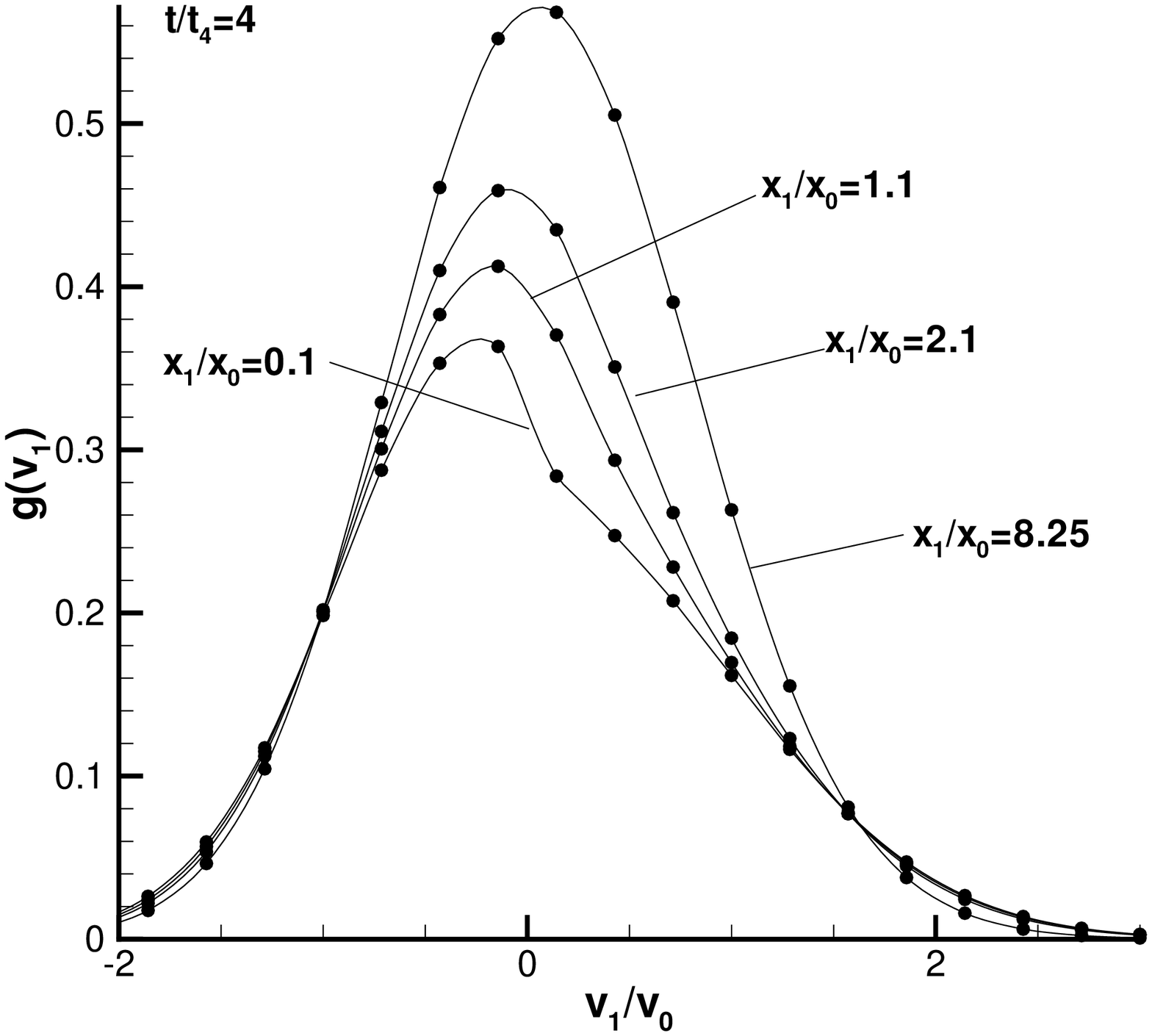}
\caption{Sudden heating: evolution of discontinuous marginal distribution near the wall.}
\end{figure}

In Table \ref{times} we show the wall time scaling results for the sudden heating example in Figure 1. This code was executed on the TACC supercomputer Lonestar, which has twelve cores per node. There is a jump in computational time when moving from a single node to two nodes, but every doubling of the node amount thereafter results in an almost exact halving of the computational time required with the previous number of nodes scowing near-perfect linear scaling for this range of nodes.
%

\begin{table}[!htbp]
\begin{tabularx}{\textwidth}{XXX}
\hline
nodes & cores & time (s)\\
\hline
1 & 12 &  203\\
2 & 24 & 235.3\\
4 & 48 &  120.8\\
8 & 96 &  61.4\\
16 & 192 &  30.9\\
32 & 384 &  15.2\\
64 & 768 &  7.7\\
\hline
\end{tabularx}
\caption{Computational time for a single timestep in sudden heating example from Figure 1.}
\centering
\label{times}
\end{table}


\section{Conclusions} \label{sec:conclusion}
We have extended the spectral method of Gamba and Tharkabhushanam to a second order scheme with a nonuniform grid in physical space, and investigated its scaling to high performance computers. The method showed nearly linear speedup when applied to a large computation problem solved on a supercomputer. However, at some point memory access will still become a problem, not with transfer between nodes but with memory access on a single node. The most expensive object to store is the six dimensional weight array $\widehat{G}(\zeta,\xi)$, and if $N$ becomes too large it may not fit on a single node's memory, significantly slowing down computation. At that point, it may be faster to simply compute the weights on the fly, as flops are much cheaper than memory accesses on the large distributed systems used in high performance computing today. Future advances in hardware such as Intel's new Many Integrated Core (MIC) nodes, which will be used in TACC's new Stampede computer starting in 2013, have many more cores and much more local memory than ever before. Future work will investigate scaling when pushing the bounds of both number of cores $p$ of the system as well as the memory requirements on each node.

\section{Acknowledgments}
This work has been supported by the NSF under grant number DMS-0636586.





\bibliographystyle{siam}   

\bibliography{Boltz}


\end{document}